\documentclass[11pt,bookmarksopen]{article} 

\usepackage[dvipdfm]{hyperref} \usepackage{amssymb} 
 \usepackage{indentfirst,latexsym, bm}
  \usepackage{amsmath}

\begin{document}

  \title{ Divergence form nonlinear nonsmooth parabolic equations with locally arbitrary growth conditions and  nonlinear maximal regularity}
\author{ Qiao-fu Zhang \\ \vspace{2mm} \footnotesize{(Academy of Mathematics and Systems Science,} \\ \footnotesize{ Chinese Academy of Sciences, Beijing 100190, P. R. China)}
     }  

\date{} \maketitle

\newcommand{\dif}{\,\mathrm{d}}

\newtheorem{dy}{ \textbf{Definition}}[section] \newtheorem{mt}{ \textbf{Proposition}}[section] \newtheorem{yl}{ \textbf{Lemma}}[section] \newtheorem{z}{ \textbf{Remark}}[section]
\newtheorem{dl}{\textbf{Theorem}}[section] \newtheorem{tl}{\textbf{Corollary}}[section] \newcommand{\pf}{\textbf{Proof\,\,\,}} \newcommand{\epf}{\hfill$\square$}

\makeatletter 
\renewcommand\theequation{\thesection%
                   .\arabic{equation}}

\abstract{
             This is a  generalization  of our prior work on the compact f\mbox{}ixed point  theory for the elliptic Rosseland-type equations.
              Inspired  by the   Rosseland equation in the conduction-radiation coupled heat transfer,
               we use the locally arbitrary growth conditions instead of the common global restricted growth conditions. Its physical meaning is: the absolute temperature should be positive and bounded.

                There exists a f\mbox{}ixed point for the  linearized map (compact and continuous in $L^2$) in a closed convex set. We also consider the nonlinear maximal regularity in Sobolev space.
}

\vspace{5mm}

\noindent\textbf{Key words: \,\,arbitrary growth conditions; f\mbox{}ixed point;}

\textbf{ Rosseland equation; nonlinear maximal regularity;}

\textbf{ nonlinear  parabolic equations; nonsmooth data. }
\section{Introduction}

 Suppose $S=(0,T)$ where  $T$ is a positive constant. Consider  the following parabolic problem: 
\begin{equation}
 \partial_t u -\mbox{div\,}[ A(u(x),x,t)\nabla u]  = 0 ,\quad \mbox{in\,}\,Q_T=\Omega\times S.
\end{equation}

The weak solution can be def\mbox{}ined as the following: f\mbox{}ind $u$,
\begin{equation}
(u-g)\in L^2(S;\,H^1_0(\Omega)),\quad (u-g)(x,0)=0 ,
\end{equation}
(so we know the boundary and initial conditions)
\begin{equation}
 \partial_t u\in L^2(S;H^{-1}(\Omega)),
\end{equation}
where $L^2(S;H^{-1}(\Omega))$ is the   dual   space of $L^2(S;\,H^1_0(\Omega))$, such that $\forall\, \varphi\in L^2(S;\,H^1_0(\Omega))$,
\begin{equation}
\langle \partial_t u,\,\varphi\rangle_{ L^2(S;\,H^1_0(\Omega))} +\iint_{Q_T} A(u(x),x)\nabla u \cdot \nabla \varphi = 0   .
\end{equation}

For the def\mbox{}initions of these spaces, see \cite{csx,intro}.

For the Rosseland equation: $A(z,x,t)=K(x,t)+z^3B(x,t)$, where $K(x,t)$ and $B(x,t)$ are symmetric and positive def\mbox{}inite.

   (1)  $K(x,t)+z^3B(x,t)$ is  positive def\mbox{}inite only in an interval for $z$.

    (2)  it doesn't satisfy the common growth and smooth conditions and there may be no $C^{2,\gamma}$ estimate (Theorem 15.11 \cite{gl}).

The problem  of the existence theory for the  Rosseland equation (also named dif\mbox{}fusion approximation) was proposed by Laitinen \cite{Laitinen} in 2002.
 It may be useful to keep this equation in mind while reading this paper.

It's a little technical to prove the existence by the f\mbox{}ixed point method in $L^\infty(Q_T)$ (or $C^0(\overline{Q}_T)$ ) \cite{amm,zqfthesis}. We will use $L^2(Q_T)$ in this paper.

F\mbox{}irstly, we make the following assumptions.

(A1)  $\Omega\subset \mathbb{R}^n$ is  a bounded Lipschitz domain. $S=(0,T)$ where  $T$ is a positive constant. $Q_T=\Omega\times S$. 


(A2) $A=(a_{ij})$. $a_{ij} =a_{ji} $.  $T_{min}\leq T_{max}$ are two constants.
\begin{equation}
   \lambda |\xi|^2 \leq a_{ij}(z,x,t)\xi_i\xi_j\leq \Lambda  |\xi|^2, \quad 0< \lambda  \leq \Lambda ,
  \end{equation}
\begin{equation}
\forall\,\,(z,x,t,\xi )\in  [T_{min},T_{max} ] \times Q_T \times \mathbb{R}^n.
\end{equation}
Here we use the Einstein   summation convention.

(A3) $\partial_pQ_T=  \{\partial \Omega \times S\}\cup \{(x,0);x\in\Omega\}$,
 \begin{equation}
g \in H^{1}( Q_T ).\quad  T_{min}\leq g(x,t)\leq T_{max},\quad \mbox{a.\,e.\,\,in}\,\,\, \partial_pQ_T.
 \end{equation}

(A4)  $A(z,x,t )$ is uniformly continuous with respect to $z$ in $\mathfrak{C}$, where
\begin{equation}
 \mathfrak{C} =\{ \varphi\in L^2(Q_T);\,\,T_{min}\leq \varphi(x,t) \leq T_{max},\,\,\mbox{a.\,e.\,\, in} \,\,Q_T\}.
\end{equation}

In other words, if $ z_i,\,z\in \mathfrak{C}$, $\|z_i-z\|_2\to 0$,
  \begin{equation}
  \sup_{1\leq p,q\leq n}\|a_{pq}(z_{i }(x,t),x,t)-a_{pq}(z(x,t),x,t)\|_2\to 0.
  \end{equation}

\begin{z}\rm
 In fact, we had considered a  general case:  parabolic equations with $(c(x)\rho(x) u)'$, nonnegative bounded mixed boundary conditions and  right-hand term $f(z,x,t)$ \cite{zqfthesis}.

If $a_{pq}$ is uniformly H\"older continuous    with respect to $z$, (A4) is natural since
  \begin{eqnarray}
  \|a_{pq}(z_{i }(x,t),x,t)-a_{pq}(z(x,t),x,t)\|_2^2&\leq& \iint_{Q_T} C |z_{i }(x,t)-z(x,t)  |^{2\alpha}
  \nonumber\\
&\leq& C\|z_i-z\|_2^2\to 0.
  \end{eqnarray}

\end{z}

\begin{dl}[\bf Parabolic spaces\rm ] Let
  \begin{equation}
W\equiv \{w\in L^2(S;\,H^1_0(\Omega));\, \partial_t w\in  L^2(S;H^{-1}(\Omega))\},
  \end{equation}
    \begin{equation}
\|w\|_W^2 =  \|w\|_{ L^2(S;\,H^1_0(\Omega))}^2 +\| \partial_t w\|_{ L^2(S;H^{-1}(\Omega))}^2,
  \end{equation}
then $($1$)$ $($page 173 \cite{csx}, page 61 \cite{intro}$)$
  \begin{equation}
W \hookrightarrow C([0,T]; L^2(\Omega )),\quad W\hookrightarrow L^2(Q_T).
  \end{equation}
The last imbedding is compact.

 $($2$)$ $($page 173 \cite{csx}$)$  $C^\infty([0,T];H^1_0 (\Omega))$ is dense in $W$.

 $($3$)$ $($Theorem 1.6 \cite{gr07}$)$ Let
 \begin{equation}
W_{c\rho}\equiv \{w\in L^2(S;\, X );\, \partial_t [c\rho(x)w]\in  L^2(S; X')\},
  \end{equation}
then   $C^\infty([0,T]; X)$ is dense in $W_{c\rho}$. For the mixed boundary conditions, we can let $X=H^1_D(\Omega)$.
\end{dl}
\begin{dl}[\bf $V^{1,0}_2(Q_T)$\rm ]
 $($page 42-44 \cite{cyz}$)$ Let
    \begin{equation}
  V^{1,0}_2(Q_T)\equiv L^2(S;\,H^1 (\Omega))\cap C([0,T];L^2( \Omega)) ,
 \end{equation}
     \begin{equation}
 \|w\|_{ V^{1,0}_2(Q_T)}^2 = \|\nabla w\|_{ L^2(Q_T;\mathbb{R}^n)}^2 + \sup_{t\in [0,T]}\| w(x,t)\|_{L^2( \Omega)}^2 ,
 \end{equation}
 then

  $($1$)$  $    H^1 (Q_T) \subset   V^{1,0}_2(Q_T)$.

  $($2$)$  If $ u(x,t) \in V^{1,0}_2(Q_T)$, $\forall\,k\in \mathbb{R}$,
    \begin{equation}
  (u-k)_+(x,t)=\max\{(u-k)(x,t),0\}\in V^{1,0}_2(Q_T) .
     \end{equation}

       $($3$)$  If $ \|u_i-u\|_{  V^{1,0}_2(Q_T)}\to 0$, then $\forall\,k\in \mathbb{R}$,
    \begin{equation}
  \|(u_i-k)_+-(u-k)_+\|_{  V^{1,0}_2(Q_T)}\to 0 .
     \end{equation}
\end{dl}

\section{Linearized map and f\mbox{}ixed point}

\begin{dl} $($Corollary 11.2 \cite{gl}$)$\label{dl:fdp} Let $\mathfrak{C}$ be a closed convex  set in a Banach space $\mathfrak{B}$ and let $\mathcal{ L}$ be a continuous mapping of $\mathfrak{C}$
into itself such that the image $\mathcal{ L}\mathfrak{C}$  is precompact. Then $\mathcal{ L}$ has a f\mbox{}ixed point.
\end{dl}
\begin{yl} The following set
\begin{equation}
 \mathfrak{C} =\{ \varphi\in L^2(Q_T);\,\,T_{min}\leq \varphi(x,t) \leq T_{max},\,\,\mbox{a.\,e.\,\, in} \,\,Q_T\}.
\end{equation}
 is a closed convex  set in the Banach space $L^2(Q_T )$.
 \end{yl}
 \pf
 Suppose $v_i\in \mathfrak{C}$, $v\in L^2(Q_T )$, $\|v_i-v\|_2\to 0$. If $v \notin \mathfrak{C}$, there exist
two constants $\delta_0>0$, $\delta_1>0$, such that the Lebesgue measure of the set $Q_0\equiv\{(x,t)\in Q_T;\, v(x,t)\geq T_{max}+\delta_0\}$
is bigger than $\delta_1>0$. Then
\begin{equation}
 \|v_i-v\|_2^2
 =\iint_{Q_T}|v_i-v | ^2\geq
 \iint_{Q_0 } |v_i-v | ^2\geq \delta_0^2 \delta_1.
\end{equation}
 It's impossible since $\|v_i-v\|_2\to 0$. Similarly, $v \geq T_{min}$ and $\mathfrak{C}$ is closed.

 \begin{equation}
   \forall\, \theta\in [0,1],\quad \theta v_1 + (1- \theta) v_2\leq \theta T_{max} + (1- \theta)T_{max}=T_{max}   .
    \end{equation}
     So $\mathfrak{C}$ is convex.
\epf

\begin{dl} If $(A1)-(A4)$ are satisf\mbox{}ied, then

$(1)$ $\forall\,  z\in \mathfrak{C} $,  there exists  a unique $w$,
\begin{equation}
(w-g)\in L^2(S;\,H^1_0(\Omega)),\quad (w-g)(x,0)=0 ,
\end{equation}
\begin{equation}
 w \in \mathfrak{C},\quad \partial_t w\in L^2(S;H^{-1}(\Omega)),
\end{equation}
such that  $\forall\,\varphi\in L^2(S;\,H^1_0(\Omega))$,
\begin{equation}
\langle \partial_t w,\,\varphi\rangle_{ L^2(S;\,H^1_0(\Omega))} +\iint_{Q_T} A(z(x,t),x,t)\nabla w \cdot \nabla \varphi  = 0   .
\end{equation}

$(2)$  Def\mbox{}ine a map $\mathcal{ L}:\, \mathfrak{C}\to  \mathfrak{C}$, $\mathcal{ L}z=w$, then $\mathcal{ L}\mathfrak{C}$ is precompact in $L^2(Q_T )$.

$(3)$  $\mathcal{ L}$ is continuous in $L^2(Q_T)$. So $\mathcal{ L}$ has a f\mbox{}ixed point in $\mathfrak{C}$.
\end{dl}
\pf
(1) For the \textit{a priori} estimate, since $H^{-1}(\Omega)\hookrightarrow L^{2}(\Omega)$ (page 55, 60 \cite{intro}),
 \begin{equation}
  \partial_t g\in L^2(Q_T)= L^2(S;L^{2}(\Omega)) \hookrightarrow  L^2(S;H^{-1}(\Omega)) ,
\end{equation}
 \begin{equation}
  v\equiv(w-g)\in W \hookrightarrow C([0,T];L^2( \Omega)),
\end{equation}
  \begin{equation}
 g \in H^1(Q_T) \hookrightarrow C([0,T];L^2( \Omega)) ,\quad w\in C([0,T];L^2( \Omega)).
\end{equation}
  \begin{equation}
  w\in L^2(S;H^{ 1}(\Omega)) ,\quad w \in V^{1,0}_2(Q_T).
\end{equation}

 Let
    \begin{equation}
 \varphi=(w-T_{max})_+\in V^{1,0}_2(Q_T)  ,
 \end{equation}
 then
 \begin{equation}
 \varphi(x,t)|_{\partial_p Q_T} =0,\quad \varphi\in L^2(S;\,H^1_0(\Omega)) .
  \end{equation}

For any  $v_i\in C^\infty([0,T];\,H^1_0(\Omega))\subset H^1(Q_T)$,   we have
\begin{eqnarray}
&&\langle \partial_t (v_i+g),\,(v_i+g-T_{max})_+\rangle_{ L^2(S;\,H^1_0(\Omega))}
 \nonumber\\
 &=&
 \iint_{Q_T} \partial_t (v_i+g)\cdot(v_i+g-T_{max})_+
 \nonumber\\
 &=&
 \iint_{Q_T} \partial_t  \frac {(v_i+g-T_{max})_+(x,t)^2}{2}
   \nonumber\\
 &=&\int_\Omega  \frac {(v_i+g-T_{max})_+(x,T)^2}{2}
    \nonumber\\
 & &
 -\int_\Omega \frac {(v_i+g-T_{max})_+(x,0)^2}{2}
 .
\end{eqnarray}

   By the density of $C^\infty([0,T];\,H^1_0(\Omega))$ in $W$,  for $v\equiv(w-g)\in W$, we can f\mbox{}ind $  \{v_i\}\subset C^\infty([0,T];\,H^1_0(\Omega))$ such that
\begin{equation}
 \|v_i-v \|_{C([0,T];L^2( \Omega)) }\leq C\|v_i-v\|_{W}\to 0,
\end{equation}
\begin{equation}
\|(v_i+g)-(v +g)\|_{ V^{1,0}_2(Q_T)}= \|v_i-v\|_{ V^{1,0}_2(Q_T)}\to 0,
\end{equation}
\begin{equation}
   \|(v_i+g-T_{max})_+-(v +g-T_{max})_+\|_{ V^{1,0}_2(Q_T)}\to 0.
\end{equation}
\begin{equation}
   \|(v_i+g-T_{max})_+-(v +g-T_{max})_+\|_{ L^2(S;\,H^1 (\Omega)) }\to 0.
\end{equation}
\begin{equation}
   v_i,\,v|_{\partial \Omega\times S }=0,\quad  g|_{\partial \Omega \times S}\leq T_{max} .
\end{equation}
\begin{equation}
   \|(v_i+g-T_{max})_+-(v +g-T_{max})_+\|_{ L^2(S;\,H^1_0 (\Omega)) }\to 0.
\end{equation}
\begin{equation}
\|\partial_t (v_i+g)- \partial_t (v+g)\|_{L^2(S;H^{-1}(\Omega)) }\to 0,
\end{equation}
\begin{eqnarray}
  & &  \int_ \Omega [(v_i+g -T_{max})_+(x,t)^2  -   (v+g -T_{max})_+(x,t)^2]
 \nonumber\\
 &=& \int_ \Omega [(v_i+g -T_{max})_+(x,t)   +   (v+g -T_{max})_+(x,t) ]
  \nonumber\\
  &&\quad
    [(v_i+g -T_{max})_+(x,t)   -   (v+g -T_{max})_+(x,t) ]
 \nonumber\\
 &\leq & \|(v_i+g -T_{max})_+(x,t) +  (v +g -T_{max})_+(x,t)\|_{L^2(\Omega)}
   \nonumber\\
  &&\quad
   \|(v_i+g -T_{max})_+(x,t)   -   (v+g -T_{max})_+(x,t)\|_{L^2(\Omega)}
  \nonumber\\
  &\leq &(\|(v_i+g -T_{max}) (x,t)\|_{L^2(\Omega)}   +  \|(v +g -T_{max}) (x,t)\|_{L^2(\Omega)} )
    \nonumber\\
  &&\quad  \|(v_i+g -T_{max}) (x,t)   -   (v+g -T_{max}) (x,t)\|_{L^2(\Omega)}
  \nonumber\\
   &\leq &(\|(v_i+g -T_{max}) (x,s)\|_{C([0,T]; L^2(\Omega))}
     \nonumber\\
  &&\quad\quad +  \|(v +g -T_{max}) (x,s)\|_{C([0,T]; L^2(\Omega))} )   \| v_i    -   v\|_{C([0,T]; L^2(\Omega))}
  \nonumber\\
  &\leq &C
     \| v_i    -   v\|_{W}\to 0 .
\end{eqnarray}
\begin{eqnarray}
&&\langle \partial_t w,\,(w-T_{max})_+\rangle_{ L^2(S;\,H^1_0(\Omega))}
 \nonumber\\
&=&\langle \partial_t (v+g),\,(v+g-T_{max})_+\rangle_{ L^2(S;\,H^1_0(\Omega))}
 \nonumber\\
 &=&\lim_{i\to \infty}\langle \partial_t (v_i+g),\,(v_i+g-T_{max})_+\rangle_{ L^2(S;\,H^1_0(\Omega))}
 \nonumber\\
  &=&  \lim_{i\to \infty}\int_ \Omega \left[\frac {(v_i+g -T_{max})_+(x,T)^2}{2} - \frac {(v_i+g-T_{max})_+(x,0)^2}{2}\right]
 \nonumber\\
 &=&  \int_ \Omega\frac {(v+g -T_{max})_+(x,T)^2}{2} -\int_ \Omega\frac {(v+g-T_{max})_+(x,0)^2}{2}
 \nonumber\\
 &=&
   \int_ \Omega\frac {(v+g -T_{max})_+(x,T)^2}{2}\geq 0 .
\end{eqnarray}
\begin{eqnarray}
 &&\iint_{Q_T} A(z(x,t),x,t)\nabla w \cdot \nabla (w-T_{max})_+
 \nonumber\\
 &=&
  \iint_{Q_T} A(z(x,t),x,t)\nabla (w-T_{max})_+ \cdot \nabla (w-T_{max})_+
  \nonumber\\
 & \geq &
  \lambda  \int_S\int_{\Omega} |\nabla (w-T_{max})_+|^2
   \nonumber\\
 & \geq &
  C(\Omega)\lambda \int_S\int_{\Omega}   (w-T_{max})_+ ^2
   .
\end{eqnarray}
\begin{eqnarray}
& &  C(\Omega)\lambda \int_S\int_{\Omega}   (w-T_{max})_+ ^2
    \nonumber\\
 & \leq &  \langle \partial_t w,\,(w-T_{max})_+ \rangle_{ L^2(S;\,H^1_0(\Omega))}
  \nonumber\\
  &&\quad\quad
  +\iint_{Q_T} A(z(x,t),x,t)\nabla w \cdot \nabla (w-T_{max})_+
   \nonumber\\
 &=& 0.
\end{eqnarray}

So $w\leq T_{max}$, a. e. in $Q_T$. Similarly, $w\in \mathfrak{C }$.

For the well-posedness (the existence, uniqueness and the estimate  in $W$)  of $w,\,\,(w-g)\in W$, we refer to Galerkin method (page 171 \cite{csx}, page 77 \cite{czm}, page 205-211 \cite{intro}; for mixed problems, see Theorem 2.2 \cite{gr07}).

(2) $\| (w-g ) \|_{W}\leq C$. $W$ can be compactly imbedded in $L^2(Q_T)$, so $\mathcal{ L}\mathfrak{C}$ is precompact in $L^2(Q_T)$.

(3)  Suppose
\begin{equation}
 z_i,\,z\in \mathfrak{C},\quad \|z_i-z\|_2\to 0,\quad \mathcal{ L}z_i=w_i , \quad \mathcal{ L}z =w.
\end{equation}

 $ W$ is a Hilbert and thus a ref\mbox{}lexive space, so there exists a subsequence $\{i_k\}$ and $v_0=(w_0-g)\in W$ such that
\begin{equation}
(w_{i_k}-g) \to (w_0-g) ,\quad \mbox{weakly\,\,in\,\,} W .
\end{equation}
\begin{equation}
 W \subset L^2(S;\,H^1_0(\Omega)),\quad (L^2(S;\,H^1_0(\Omega)))'\subset W' .
\end{equation}
\begin{equation}
(w_{i_k}-g) \to (w_0-g) ,\quad \mbox{weakly\,\,in\,\,}  L^2(S;\,H^1_0(\Omega)) .
\end{equation}
\begin{equation}
\nabla (w_{i_k}-g) \to\nabla (w_0-g) ,\quad \mbox{weakly\,\,in\,\,}  L^2(Q_T;\mathbb{R}^n ) .
\end{equation}
\begin{equation}
\nabla  w_{i_k}  \to\nabla  w_0  ,\quad \mbox{weakly\,\,in\,\,}  L^2(Q_T;\mathbb{R}^n ) .
\end{equation}
\begin{equation}
  \|w_{i_k} -g -w_0+g\|_2\to 0 ,\quad \|w_{i_k} -w_0\|_2\to 0  .
\end{equation}
 \begin{equation}
\partial_t (w_{i_k}-g) \to \partial_t (w_0-g) ,\quad \mbox{weakly\,\,in\,\,}  L^2(S;\,H^{-1}(\Omega)) .
\end{equation}
 \begin{equation}
\partial_t  g\in  L^2(Q_T)\subset  L^2(S;\,H^{-1}(\Omega)) .
\end{equation}
 \begin{equation}
\partial_t  w_{i_k}  \to \partial_t  w_0 ,\quad \mbox{weakly\,\,in\,\,}  L^2(S;\,H^{-1}(\Omega)) .
\end{equation}

 $\forall\, \phi\in C^{\infty} ([0,T];\,C^{\infty}_0(\Omega))$, using the natural map into  its second dual (page 89 \cite{Conway}),  
  \begin{eqnarray}
&&\langle F(\phi), \partial_t  w_{i_k} -\partial_t  w_0 \rangle_{L^2(S;\,H^{-1} (\Omega))    }
 \nonumber\\ & \equiv
 &
\langle \partial_t  w_{i_k} -\partial_t  w_0 , \phi\rangle_{L^2(S;\,H^{ 1}_0 (\Omega))    }
.
\end{eqnarray}
  \begin{eqnarray}
&&
\langle F(\phi), \partial_t  w_{i_k} -\partial_t  w_0 \rangle_{L^2(S;\,H^{-1} (\Omega))    }\to 0,
 \nonumber\\ &  \Rightarrow
 &
\langle \partial_t  w_{i_k} -\partial_t  w_0 , \phi\rangle_{L^2(S;\,H^{ 1}_0 (\Omega))    }\to 0
.
\end{eqnarray}

   \begin{eqnarray}
 &&|\iint_{Q_T} [A(z_{i_k} (x,t),x ,t)\nabla w_{i_k} -A(z  (x,t),x,t )\nabla w_0] \cdot\nabla \phi |
\nonumber\\ &\leq& |\iint_{Q_T} [A(z_{i_k} (x,t),x,t )\nabla w_{i_k} -A(z  (x,t),x,t )\nabla w_{i_k}] \cdot\nabla \phi | 
\nonumber\\&&\,+\,|\iint_{Q_T} [A(z  (x,t),x,t )\nabla w_{i_k} -A(z  (x,t),x ,t)\nabla w_0]
\cdot\nabla \phi | 
\nonumber\\ &=&
 |\iint_{Q_T}[A(z_{i_k} (x,t),x,t )  -A(z  (x,t),x,t )] \nabla w_{i_k}\cdot\nabla \phi | 
 \nonumber\\&&
 \,+\,|\iint_{Q_T}[\nabla w_{i_k} - \nabla w_0] \cdot A(z  (x,t),x,t
)^\top\nabla \phi | 
\nonumber\\ &\leq& C \sup_{1\leq p,q\leq n}\|a_{pq}(z_{i_k}(x,t),x,t)-a_{pq}(z(x,t),x,t)\|_2+\epsilon(i_k) 
\nonumber\\ &\to &0
 .
\end{eqnarray}
  \begin{eqnarray}
&&  
\langle \partial_t  w_{i_k}  , \phi\rangle_{L^2(S;\,H^{ 1}_0 (\Omega))    }
\nonumber\\&&\,+\,\iint_{Q_T}  A(z_{i_k} (x,t),x,t )\nabla w_{i_k}\cdot\nabla \phi=0
.
\end{eqnarray}
  \begin{eqnarray}
&&
\langle \partial_t  w_{0}  , \phi\rangle_{L^2(S;\,H^{ 1}_0 (\Omega))    }
\nonumber\\&&\,+\,\iint_{Q_T} A(z  (x,t),x,t )\nabla w_{0}\cdot\nabla \phi=0
.
\end{eqnarray}

From the density of $C^{\infty} ([0,T];\,C^{\infty}_0(\Omega))$ in $L^2(S;\,H^{ 1}_0 (\Omega))$, $\forall\,\,\varphi\in L^2(S;\,H^{ 1}_0 (\Omega))$,
  \begin{eqnarray}
&&
\langle \partial_t  w_0 , \varphi\rangle_{L^2(S;\,H^{ 1}_0 (\Omega))    }
\nonumber\\&&\,+\,\iint_{Q_T} A(z  (x,t),x,t )\nabla w_{0}\cdot\nabla \varphi=0
.
\end{eqnarray}

For the boundary condition,
\begin{equation}
 (w_0-g)\in L^2(S;\,H^1_0(\Omega)) ,\quad (w_0-g)|_{\partial \Omega \times S}=0.
\end{equation}

For the initial condition,
 $\forall \,\,\psi(x)\in L^2(\Omega) $, we can def\mbox{}ine a linear functional on $W$, 
\begin{equation}
 \langle \Psi  , h\rangle_{ W  } \equiv \int_ \Omega h(x,0)\psi,\quad \forall \,\, h(x,t)\in W.
\end{equation}

This functional is bounded since
  \begin{eqnarray}
|\langle \Psi  , h\rangle_{ W  }|&=&
|\int_ \Omega h(x,0)\psi|
\nonumber\\
&\leq & \|h(x,0)\|_2\|\psi\|_2\leq \|\psi\|_2\sup_{s\in [0,T]}\|h(x,s)\|_2
\nonumber\\
&= & \|\psi\|_2  \|h(x,t)\|_{C([0,T];\, L^2(\Omega))}
\nonumber\\
&\leq & C\|h(x,t)\|_{W}
.
\end{eqnarray}

Since
\begin{equation}
(w_{i_k}-g) \to (w_0-g) ,\quad \mbox{weakly\,\,in\,\,} W ,
\end{equation}
$\forall \,\, \psi(x)\in L^2(\Omega)$,
\begin{equation}
 \langle \Psi  , (w_{i_k}-g) - (w_0-g)\rangle_{ W  } \equiv \int_ \Omega [(w_{i_k}-g) - (w_0-g)](x,0)\psi\to 0  .
\end{equation}

From the Riesz Representation Theorem in $ L^2(\Omega)$,  
\begin{equation}
   (w_{i_k}-g)(x,0) \to (w_0-g) (x,0) ,\quad \mbox{weakly\,\,in\,\,}L^2(\Omega) .
\end{equation}

Note that from the initial condition,
\begin{equation}
 ( w_{i_k}-g)(x,0) =0,\quad \mbox{in}\,\, L^2(\Omega) .
\end{equation}
\begin{equation}
   (w_{i_k}-g)(x,0) \to 0 ,\quad \mbox{strongly\,\,in\,\,}L^2(\Omega) .
\end{equation}
\begin{equation}
 ( w_{0}-g)(x,0) =0,\quad \mbox{in}\,\, L^2(\Omega) .
\end{equation}

To sum up, $w_0=\mathcal{ L}z$: $w_0$ satif\mbox{}ies the linearized equation and the initial-boundary conditions. 

Since the solution is unique from the step (1), $w_0=\mathcal{ L}z =w$. So  $\|w_{i_k} -w \|_2\to 0 $.
 Each subsequence of $  \{\|w_{i } -w \|_2 \}$ has a sub-subsequence which  converges to $0$, so $\|w_{i } -w \|_2\to 0 $. We have obtain the continuity of $\mathcal{ L}$.

  From Theorem \ref{dl:fdp}, there exists a f\mbox{}ixed point.
\epf 
\begin{z}
\rm 
For the continuity of $\mathcal{ L}$ in $C^0(\overline{Q}_T)$, we can use the well-known De Giorgi-Nash estimate: $\{w_i\}$ is bounded in $C^{2\alpha ,\alpha} (\overline{Q}_T)$ if
$g\in C^{2\alpha_0,\alpha_0} ( \partial_p Q_T )$ and $\Omega$ is an (A) domain (page 145 \cite{cyz}).

Then from the Arzel$\grave{\rm{a}}$-Ascoli Lemma, $\|w_{i_k} -w _0\|_{C^0(\overline{Q}_T )}\to 0$. By the same method, $w_0=w$ and $\|w_{i } -w  \|_{C^0(\overline{Q}_T)}\to 0$.

From the linear maximal regularity \cite{gri,gr07},  a natural conjecture is:    $\mathcal{ L}$ is continuous in $C^{2\alpha,\alpha} (\overline{Q}_T)$ and $W$. 
\end{z}

\section{Nonlinear maximal regularity}

For the linear  parabolic/ellptic equations with nonsmooth data, the theory of maximal regularity has been established \cite{gri,gr07,grell,gr}. In brief, maximal regularity is about the smoothness of
the data-to-solution-map \cite{gr}. This smooth dependence has its physical meaning: many physical processes  are stable with respect to the parameters (except the chaos and critical theory).
 For the mathematicians, "the door is open to apply the powerful theorems of dif\mbox{}ferential calculus"(\cite{gr}, e.g.  the Implicit Function Theorem).

In the following, we will discuss the continuous dependence (between the solutions and the data) for the    parabolic equations with locally arbitrary growth conditions (e.g. Rosseland-type).

\section{Acknowledge }

This work is supported by the National Nature Science Foundation of China (No. 90916027). This is a part  of my PhD thesis \cite{zqfthesis} in AMSS,  Chinese Academy of Sciences, and a
simplif\mbox{}ication of our prior paper \cite{amm}.  So I will  thank my advisor Professor Jun-zhi Cui (he is also a member of the Chinese Academy of Engineering) and the referees for their careful
reading and helpful comments. My E-mail is:  zhangqf@lsec.cc.ac.cn.

\end{document}